\documentclass[11pt,a4paper]{article}
\usepackage{amsmath,amsthm,amsfonts,bbm}
\usepackage [curve,arrow,frame]{xy}
\usepackage{graphicx}

\newtheorem{thm}{Theorem}[section]

\newtheorem{lem}[thm]{Lemma}

\theoremstyle{definition}
\newtheorem{defn}[thm]{Definition}
\newtheorem{rem}[thm]{Remark}

\newcommand{\lemref}[1]{Lemma~\ref{#1}}

\newcommand {\R}{\mathbbm R}
\newcommand {\Z}{\mathbbm Z}
\newcommand {\N}{\mathbbm N}

\newcommand {\F}{\mathcal F}
\newcommand {\U}{\mathcal U}

\newcommand {\g} {\mathfrak g}

\begin{document}
\title{Quantization of $SL(2,\R)^*$ as Bialgebra}
\date{November 22, 2002} 
\author{Markus R. Engeli}
\maketitle

\abstract{We quantize the Poisson-Lie group $SL(2,\R)^*$ 
as a bialgebra using the product of Kontsevich. 
The coproduct is a deformation of the coproduct that comes
from the group structure.
The resulting bialgebra structure is isomorphic to the quantum universal enveloping algebra $U_hsl(2,\R)$.}

\section{Idea and Results}
\subsection{Idea}
Let $G$ be a semi-simple Lie group and $\g$ its Lie algebra.
If one chooses a Cartan decomposition of $\g$, there is a natural $r$-matrix in $\g\otimes\g$ which defines a Lie bialgebra structure on $\g$ resp. its dual $\g^*$: $(\g,[\,],\delta)$ resp. $(\g^*,[\,]^*,\delta^*)$.
\begin{center}
$
\begin{array}[3]{ccc}
(S(\g),\Delta) & = & (\F(\g^*),\cdot,\Delta)\\[1ex]
\downarrow & & \downarrow\\[1ex]
(\U(\g)[[\varepsilon]],\Delta_\varepsilon) & & (\F(\g^*),\cdot,\Delta_h)\\[1ex]
\downarrow & & \downarrow\\[1ex]
(\U_{h,\varepsilon}(\g),\Delta_{h,\varepsilon}) & \overset{?}{=} & (\F(\g^*),*_\varepsilon,\Delta_{h,\varepsilon})
\end{array}
$
\end{center}
\vspace{2mm}
The diagram shows 6 bialgebras. The bialgebras on the LHS are well known from quantum group theory: $S(\g)$ is the free commutative algebra generated by $\g$. Its coproduct is given by $\Delta(x)=1\otimes x + x\otimes 1$ for $x\in\g$. $\U(\g)[[\varepsilon]]$ is the universal enveloping algebra of $\g$ but with Lie bracket multiplied by $\varepsilon$. For elements in $\g$ the coproduct is the same as for $S(\g)$ and the coproduct is extended as an algebra homomorphism. The third bialgebra on the LHS, $\U_{h,\varepsilon}(\g)$, is the quantum universal enveloping algebra of $\g$ which is a deformation of $\U(\g)[[\varepsilon]]$. The structure of this algebra will be given in the next subsection.
The bialgebras on the RHS are algebras of functions on $\g^*$. In our example we will consider an algebra of finite products of some generating functions.

The two equivalent bialgebras on the top are commutative and cocommutative. The algebra $S(\g)$ can also be understood as the algebra of polynomial functions on $\g^*$. This explains why the bialgebras on the top are equivalent.\par
The bialgebras in the middle and the bottom of the diagram are deformations of the bialgebras on the top. The ``direction'' of the deformation is given by the Lie bialgebra structure of $\g$ resp. $\g^*$.\par
On the LHS we first deform the product to a non-commutative one (parameter $\varepsilon$) and in the second step also the coproduct (parameter $h$). All this is well known from quantum group theory.\par
On the RHS we first deform the coproduct and the Poisson structure on $\g^*$ ($\delta^*$ defines a Poisson bracket on $\F(\g^*)$) so that the algebra can be identified via the exponential map with the algebra of functions on the Poisson-Lie group $(G^*,\{\}_h)$. In the second step the product is deformed in the direction of the Poisson structure using the formula of Kontsevich. In this quantization the coproduct is held fixed on the generators of $\F(\g^*)$ \footnote{In ``our'' case the coproduct is also not deformed for elements in $\g\subset\F(\g^*)$.} but as it has to be compatible with $*_\varepsilon$ it gets also deformed.\par
The question which now naturally arises is wheter the two deformed algebras on the bottom of the diagram are isomorphic or not.

\subsection{Results}
We worked out explicitly the case $\g=sl(2,\R)$, and we showed that the resulting deformed bialgebras are isomorphic.
The algebra is generated by elements $\xi_1$, $\xi_2$ and $\xi_3$ satisfying the following relations:
\begin{align*}
  [\xi_1,\xi_2] &= \phantom{-}2\varepsilon \xi_2\\
  [\xi_1,\xi_3] &= -2\varepsilon \xi_3\\
  [\xi_2,\xi_3] &= \varepsilon\frac{\sinh(h\,\xi_1)}{\sinh(h)}
\end{align*}
and
\begin{align*}
  \Delta \xi_1 &= 1\otimes \xi_1 + \xi_1\otimes 1\\
  \Delta \xi_2 &= \xi_2\otimes e^{-h \xi_1/2}+e^{h \xi_1/2}\otimes \xi_2\\
  \Delta \xi_3 &= \xi_3\otimes e^{-h \xi_1/2}+e^{h \xi_1/2}\otimes \xi_3
\end{align*}

If we take the limit $h\to 0$ we get the bialgebra in the middle of the LHS of the diagram. The limit $\varepsilon\to 0$ corresponds to the bialgebra in the middle of the RHS of the diagram. And of course taking both limits gives the bialgebra on the top of the diagram.
\\\par
{\bf Acknowledgements:} I am greatly indebted to my advisor G. Felder who made it possible to me to write this article. I am also grateful to A.S. Cattaneo for very helpful discussions.

\section{Lie Bialgebra}
A basis for $sl(2,\R)$ is given by the matrices
\begin{equation*}
  H=
  \begin{pmatrix}
    1 & 0\\
    0 & -1
  \end{pmatrix}
  \qquad
  X^+=
  \begin{pmatrix}
    0 & 1\\
    0 & 0
  \end{pmatrix}
  \qquad
  X^-=
  \begin{pmatrix}
    0 & 0\\
    1 & 0
  \end{pmatrix}
\end{equation*}
The Lie algebra structure is characterized by
\begin{equation*}
  [H,X^\pm] = \pm 2 X^\pm \qquad [X^+,X^-] = H
\end{equation*}
The cobracket (or cocommutator) is the coboundary of the r-matrix
\begin{equation*}
  r = \;X^+\wedge X^- = \,(X^+\otimes X^- - X^-\otimes X^+)
\end{equation*}
This gives the standard Lie bialgebra structure on $sl(2,\R)$:
\begin{equation*}
  \delta(H)=0 \quad\text{and}\quad \delta(X^\pm)=X^\pm\wedge H
\end{equation*}

\section{Poisson-Lie Group}
First we explain the notions of Poisson-Lie group and Lie bialgebra. A more detailed introduction can be found in \cite{DF} or \cite{CP}.
\begin{defn}
  A {\em Poisson-Lie group} is a Lie group with a Poisson structure that is compatible with the Lie structure in the sense that the group multiplication ($m: G\times G\to G$) is a Poisson map.
\end{defn}
The infinitesimal neighborhood of the identity element of a Poisson-Lie group gives us the structure of a Lie bialgebra.
\begin{defn}
A {\em Lie bialgebra} is a Lie algebra with the additional structure of a linear map $\delta:\g\to\g\otimes\g$, called the cobracket (or cocommutator). The dual of the cobracket $\delta^*$ has to be a Lie bracket on $\g^*$ and the cobracket $\delta$ has to be a 1-cocycle in the Lie algebra cohomology of $\g$ with values in $\g\otimes\g$ with adjoint action. 
\end{defn}
In the case of Poisson-Lie groups, the cobracket can be given as follows. Let $\alpha$ be the Poisson bivector field on $G$. Then define $w:G\to\g\otimes\g$ as $w(g):=R_g^*\alpha(g)$ where $R_g:G\to G$ is the right translation by the element $g\in G$. Then $\delta$ is given as the derivative of $w$ at the identity, i.e. $\delta = dw(\mathbbm{1}) : \g\to\g\otimes\g$.

\begin{rem} Let $\g$ be a Lie bialgebra
  \begin{enumerate}
  \item Also $\g^*$ is a Lie bialgebra: its bracket is the dual of the cobracket of $\g$ and its cobracket is the dual of the bracket of $\g$.
  \item Each Lie bialgebra can be integrated to a Poisson-Lie group. This means that there is a Poisson-Lie group, so that the Lie bialgebra describes its infinitesimal structure at the identity.
  \end{enumerate}
\end{rem}

We take the dual Lie bialgebra\footnote{When defining duality we could scale the dual bracket and dual cobracket by a constant factor. This will not be interesting for the calculations that follow but we will use that fact in the last sections for the interpretation of our results.} and integrate it to the following Lie group (see \cite{CP}, example 2.2.10)
\begin{align*}
  SL(2,\R)^* &= \left\{ \left. \left[ 
        \begin{pmatrix} a^{-1} & 0\\ b & a \end{pmatrix}, 
        \begin{pmatrix} a & c\\ 0 & a^{-1} \end{pmatrix} \right] 
    \ \right| \ a\in\R_+, b,c \in\R
  \right\} \\[2ex]
&\subset SL(2,\R)\times SL(2,\R)
\end{align*}
\begin{lem}
  Let G be a Lie group and $\g$ its Lie algebra. Assume that $\exp:\g\to G$ is a bijective map and let $\delta$ be a cobracket on $\g$. Then the integration of the cobracket leads to the Poisson structure
  \begin{equation*}
    \alpha(e^X)=(R_{e^X})_*\int_0^1(Ad_{e^{sX}}\otimes Ad_{e^{sX}})\delta(X)ds
  \end{equation*}
  which makes $G$ to a Poisson-Lie group. $R$ denotes the right translation on the group.
\end{lem}
\begin{proof}
  Let $G$ be an arbitrary Poisson-Lie group with Poisson bivector field $\alpha$ and Lie algebra $\mathfrak{g}$. The pullback of the Poisson bivector field
  \begin{equation*}
    w(g):=R_g^*\alpha(g)\in\mathfrak{g}\otimes\mathfrak{g}
  \end{equation*}
  satisfies the Poisson-Lie condition
  \begin{equation*}
    w(gh) = (Ad_g\otimes Ad_g)w(h)+w(g)
  \end{equation*}
  with $g,h\in G$. In particular this means that $w(\mathbbm{1})$ is always zero. We now set $g=e^{sX}$ and $h=e^{tX}$ ($X\in\mathfrak{g}$) and take the derivative of the above condition with respect to $t$ in $t=0$:
  \begin{equation*}
    \frac{d}{ds}w(g(s)) = (Ad_{g(s)}\otimes Ad_{g(s)})\delta(X)
  \end{equation*}
  Integration gives:
  \begin{equation*}
    w(e^X)=\int_0^1(Ad_{e^{sX}}\otimes Ad_{e^{sX}})\delta(X)ds
  \end{equation*}
  This formula is well defined if the exponential map is bijective.
\end{proof}
The integration of the cobracket leads to the following
Poisson bivector field on $SL(2,\R)$:
\begin{equation*}
  \alpha(a,b,c) = -\frac{1}{2} ab\;\partial_a\wedge\partial_b + \frac{1}{2} ac\;\partial_a\wedge\partial_c + \frac{a^4-1}{a^2}\;\partial_b\wedge\partial_c
\end{equation*}
In the coordinates $x_1:=\ln a$, $x_2:=-b$ and $x_3:=c$ the Poisson
bivector field takes the form
\begin{equation*}\label{EQ:alpha}\tag{$*$}
  \alpha(x_1,x_2,x_3) = \underbrace{x_2\,\partial_1\wedge\partial_2}_{=:\alpha_1}  \underbrace{- x_3\,\partial_1\wedge\partial_3}_{=:\alpha_2} + \underbrace{4\sinh(2x_1)\,\partial_2\wedge\partial_3}_{=:\alpha_3}
\end{equation*}

\section{Kontsevich's Product and Gauge Transformation}
Kontsevich gave a formula for the deformation quantization of a Poisson manifold (\cite{K}). His product looks as follows:
\begin{equation*}
  f*g = \sum_{n=0}^\infty\varepsilon^n\sum_{\Gamma\in G_n}\omega_\Gamma B_{\Gamma,\alpha}(f,g)
\end{equation*}
where $G_n$ is the set of all diagrams of a certain type (admissible dia\-grams) with $n+2$ vertices. The diagrams represent bidifferential operators. Each vertex in a diagram stands for the Poisson bivector field multiplied by $\varepsilon$.

We now consider the algebra $\F$ of functions on $SL(2,\R)^*$ generated by the elements $x_1$, $x_2$, $x_3$ and $e^{\pm x_1}$.
Kontsevich's product on these generators can be written in the following form:
\begin{align*}\label{EQ:gr1}
  x_1*x_1 &= x_1^2 + c\varepsilon^2 & x_1*x_2 &= x_1x_2 + \varepsilon x_2\\\tag{a}
  x_2*x_2 &= x_2^2                  & x_1*x_3 &= x_1x_3 - \varepsilon x_3\\
  x_3*x_3 &= x_3^2&&
\end{align*}
\begin{equation*}\label{EQ:gr2}\tag{b}
  x_2*x_3 = x_2x_3 + \varepsilon A(\varepsilon^2)\sinh(2x_1) + \varepsilon^2 B(\varepsilon^2)\cosh(2x_1)
\end{equation*}
\vspace{0mm}
\begin{equation*}\label{EQ:gr3}\tag{c}
  \begin{aligned}
    e^{\pm x_1}*x_2 &= \sum_{k=0}^\infty \frac{(\pm 2\varepsilon)^k}{k!}\widehat{B}_k \, x_2 \, e^{\pm x_1}\\
    e^{\pm x_1}*x_3 &= \sum_{k=0}^\infty \frac{(\mp 2\varepsilon)^k}{k!}\widehat{B}_k \, x_2 \, e^{\pm x_1}
  \end{aligned}
\end{equation*}
where $c\in\R\setminus\{0\}$, $A(\varepsilon^2)$ and $B(\varepsilon^2)$ are formal power series in $\varepsilon^2$ and $(-1)^k\widehat{B}_k=:B_k$ are the Bernoulli numbers that are defined by $\sum_{k=0}^\infty B_k\, x^k/k! = x/(\exp(x)-1)$. We omitted the products $x_1*e^{\pm x_1}$ and $e^{x_1}*e^{-x_1}$. We will later consider a gauge equivalent product (see \lemref{LEM:gauge}) for which these products are trivial.
To compute these products, it is convenient to consider 3 types of vertices that correspond to the 3 terms $\alpha_1$, $\alpha_2$ and $\alpha_3$ ($\alpha=\alpha_1+\alpha_2+\alpha_3$, see equation (\ref{EQ:alpha})) of the Poisson bivector field.\par
If we consider $*$-products of terms that are polynomial in $x_1$, $x_2$ and $x_3$, we see that the $\alpha_1$- and $\alpha_2$-vertices reduce the power with respect to $x_1$ by one but preserve the powers with respect to $x_2$ and $x_3$. On the other hand the $\alpha_3$-vertex reduces the power of $x_2$ and $x_3$ by one but contributes with a power series with respect to $x_1$.

By counting powers one readily sees that the first five products (\ref{EQ:gr1}) are at most of second order in $\varepsilon$. Further one verifies that for the product $x_2*x_3$ (\ref{EQ:gr2}) only diagrams with exactly one $\alpha_3$-vertex give non-vanishing contributions: diagrams with no $\alpha_3$-vertex can be split up in a part that acts only on $x_2$ and a part that acts only on $x_3$. Such diagrams are of weight 0. By counting the powers with respect to $x_2$ or $x_3$, we see that more than one $\alpha_3$-vertex is not possible. Therefore, only the following types of diagrams are of interest:

\begin{center}
  \includegraphics[bb = 130 430 415 715]{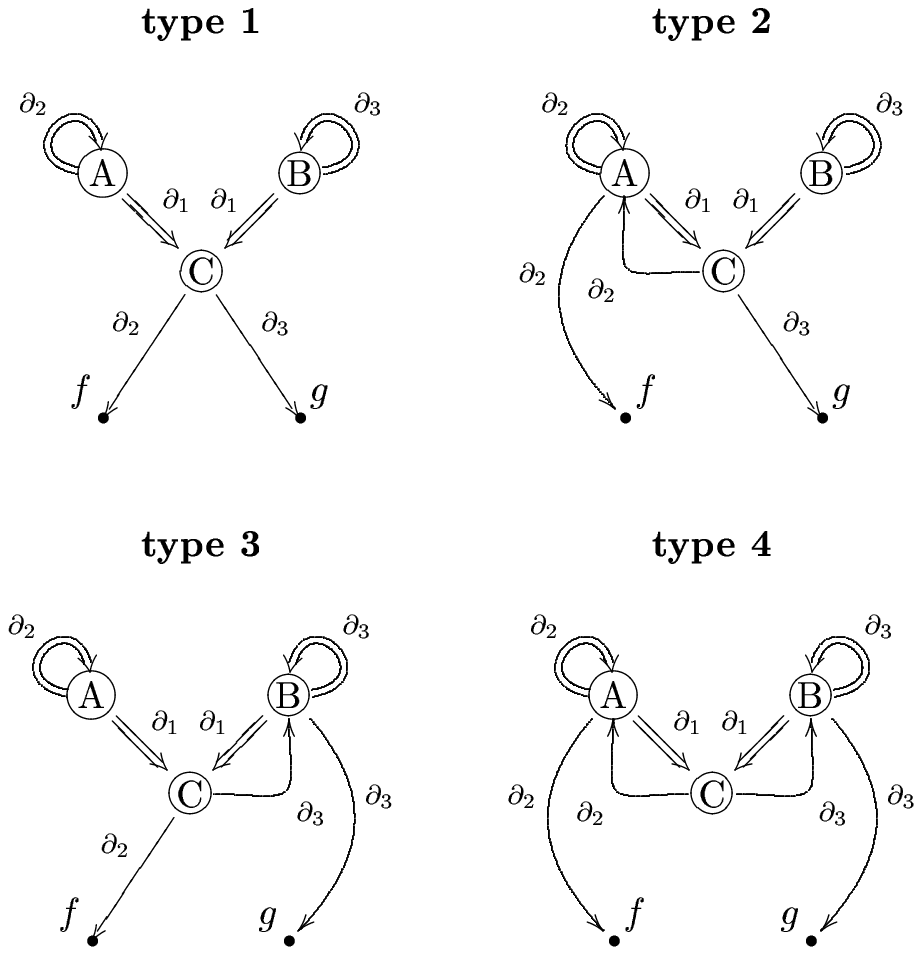}
\end{center}

C stands for the $\alpha_3$-vertex and A resp. B stand for many $\alpha_1$ resp. $\alpha_2$ vertices. Each normal arrow represents a single derivative and the double arrows stand for many derivatives (it is easy to figure out {\em how}\, many if the number of vertices is given). Since all derivatives with respect to $x_1$ act on the $\alpha_3$-vertex, all diagrams with odd (even) order in $\varepsilon$ are proportional to $\sinh(2x_1)$ ($\cosh(2x_1)$).\\

Now we consider the last two products (\ref{EQ:gr3}). For the product $e^{x_1}*x_2$ we have to compute diagrams of the following type:

\begin{center}
  \includegraphics[bb = 200 590 395 695]{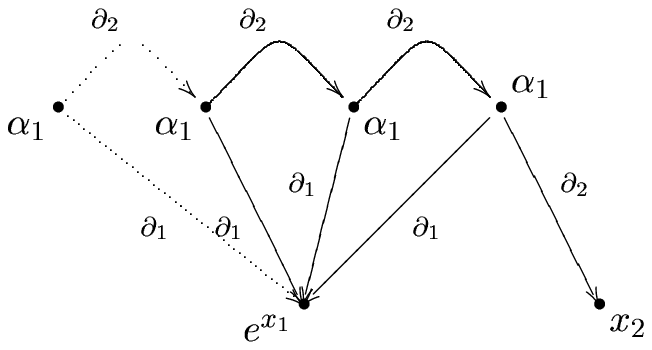}
\end{center}

The weights of diagrams of this type have already been computed by V. Kathotia \cite{Ka}. His result is that the weights are essentialy given by the Bernoulli numbers.\\

We will see in the next section that it would be convenient
to have $x_1^{*n}=x_1^n$ for all $n\in\N$. Moreover this spares 
us the computation of the products $x_1*e^{\pm x_1}$ and $e^{x_1}*e^{-x_1}$.
The solution is given in the following lemma:

\begin{lem}\label{LEM:gauge}
There is a differential operator of the form
\begin{equation*}
  D=\sum\limits_{m\in 2\N} a_m\varepsilon^m\partial_1^m
\end{equation*}
with $a_m\in\R\ \forall m\in2\N$ $(\text{and } a_0=1)$, so that the new product defined as
\begin{equation*}
  f*_{new}g = D^{-1}(Df*_{old}Dg)
\end{equation*}
obeys $x_1^{*n}=x_1^n$. $(x_1^{*n}:=\underbrace{x_1*\dots*x_1}_{\text{n factors}})$
\end{lem}
\begin{proof}
  For convenience we will write $x$ instead of $x_1$ in this proof.
  In the first step we consider the product $x^{n-1}*x$. By counting powers with respect to $x_2$ or $x_3$, one sees that no diagrams with $\alpha_3$-vertices can appear in this product. It is also not hard to see that only diagrams with only $\alpha_1$-vertices and diagrams with only $\alpha_2$-vertices give non-zero contributions. Because the $\alpha_1$- and the $\alpha_2$-vertices have opposite signs, the diagrams with an odd number of vertices cancel out. Therefore, we can write the product as 
  \begin{equation}\label{Eq_bk}
    x^{n-1}*x = x^n +
    \sum_{2\le k\le n\atop k\ even}^\infty
    \varepsilon^k\frac{(n-1)!}{(n-k)!}\; b_k \,x^{n-k}
  \end{equation}
  By the preceding arguments, it is clear that the multiple $*$-product $x^{*n}$ only consists of even orders in $\varepsilon$:
  \begin{equation}\label{Eq_cnk}
    x^{*n} = \sum_{0\le k\le n\atop k\ even} \varepsilon^k c^n_k\,x^{n-k}
    \qquad\text{with}\ c^n_0=1\ \forall n\in\N
  \end{equation}
  If we know all the $b_k$, the $c^n_k$ can be found recursively.
  We compute $x^{*(n-1)}*x$ using first equation (\ref{Eq_cnk}) and then equation (\ref{Eq_bk}):
  \begin{align*}
    x^{*n} &= \sum_{0\le l\le n-1\atop l\ even}
    \varepsilon^l c^{n-1}_l x^{n-l-1}*x\\
    &= \sum_{0\le l\le n-1\atop l\ even}
    \left( \varepsilon^lc^{n-1}_lx^{n-l} + 
      \sum_{2\le m\le n-l\atop m\ even} \varepsilon^{l+m}
      c^{n-1}_l \frac{(n-l-1)!}{(n-l-m)!} b_m x^{n-l-m}
    \right)
  \end{align*}
  If we compare this expression with equation (\ref{Eq_cnk}), we find:\\[2ex]
  $k=0$:
  \begin{equation*}
    c^n_0=c^{n-1}_0=1
  \end{equation*}
  $k>0$:
  \begin{equation*}
    c^n_k = c^{n-1}_k + \sum_{0\le l\le k-2 \atop l\ even} 
    c^{n-1}_l \frac{(n-l-1)!}{(n-k)!} b_{k-l}
  \end{equation*}
  On the other hand we want $D$ to satisfy
  \begin{equation*}
    Dx^n = \sum_{k\ge 0\atop k\ even} \varepsilon^k a_k \frac{n!}{(n-k)!} x^{n-k} = x^{*n}
  \end{equation*}
  Because of equation (\ref{Eq_cnk}), this leads to
  \begin{equation*}
    c^n_k = a_k\frac{n!}{(n-k)!}
  \end{equation*}
  If we use this result in the recursive formula for $c^n_k$, we easily find:
  \begin{equation*}
    a_k = \frac{1}{k}\sum_{2\le m\le k\atop m\ even}a_{k-m}b_m
  \end{equation*}
  Fortunately, the equation does not depend on $n$ any more.
\end{proof}

We now simply write $*$ for the {\em new} product. The products of our generators then look as follows:
\begin{align*}
  x_j*x_j &= x_j^2 & e^{\pm x_1}*e^{\pm x_1} &= e^{\pm 2x_1}\\
  x_1*x_2 &= x_1x_2 + \varepsilon x_2 & e^{\pm x_1}*e^{\mp x_1} &= 1 \\
  x_1*x_3 &= x_1x_3 - \varepsilon x_3 & x_1*e^{\pm x_1} &= x_1 e^{\pm x_1}
\end{align*}
\vspace{1mm}
\begin{equation*}
  x_2*x_3 = x_2x_3 + \varepsilon \widetilde{A}(\varepsilon^2)\sinh(2x_1) + \varepsilon^2 \widetilde{B}(\varepsilon^2)\cosh(2x_1)
\end{equation*}
\vspace{1mm}
\begin{equation}\label{EQ:cepsilon}
  \begin{aligned}
    e^{\pm x_1}*x_2 &= c(\pm\varepsilon)\, x_2 \, e^{\pm x_1}\\
    e^{\pm x_1}*x_3 &= c(\mp\varepsilon)\, x_2 \, e^{\pm x_1}
  \end{aligned}
\end{equation}
\\[1ex]
where $\textstyle\frac{\widetilde{A}(\varepsilon^2)}{A(\varepsilon^2)}=\frac{\widetilde{B}(\varepsilon^2)}{B(\varepsilon^2)}=\sum\limits_{k\ge 0\atop even}(2\varepsilon)^ka_k$ and $c(\pm\varepsilon):=\sum_{k=0}^\infty \frac{(\pm 2\varepsilon)^k}{k!}\,\widehat{B}_k$. The last two products do not change under the gauge transformation because the vertices that appear in their diagrams do not depend on $x_1$: $D^{-1}(De^{\pm x_1}*Dx_j) = D^{-1}\circ D (e^{\pm x_1}*Dx_j)=e^{\pm x_1}*x_j$ for $j=2,3$.

To compute the opposite products such as $x_3*x_2$ and $x_2*e^{\pm x_1}$ we use the fact that $f*g=g*f|_{\varepsilon\mapsto -\varepsilon}$ for Kontsevich's product\footnote{The gauge transformation preserves this identity as the differentialoperator $D$ is even in $\varepsilon$.}. In particular this means that the commutator of two functions is equal to 2 times the odd (with respect to $\varepsilon$) terms of their product. We finally give a useful relation for the products in (\ref{EQ:cepsilon}) and their opposite products:
\begin{equation}\label{EQ:relation}
  \begin{aligned}
    e^{\pm x_1}*x_2 &= e^{\pm 2\varepsilon}x_2*e^{\pm x_1}\\
    e^{\pm x_1}*x_3 &= e^{\mp 2\varepsilon}x_3*e^{\pm x_1}
  \end{aligned}
\end{equation}
These relations are a direct consequence of the following relation satisfied by the Bernoulli numbers
\begin{equation*}
  \sum_{k=0}^n (-1)^k\binom{n}{k}\widehat{B}_k = \widehat{B}_n
\end{equation*}

\section{Coproduct}

We consider the algebra $\widetilde\F = \R\{x_1,x_2,x_3,e^{\pm x_1}\}[[\varepsilon]]$ of non-commutative words in the generators. We define the coproduct for this algebra $\widetilde\Delta : \widetilde\F \to \widetilde\F\otimes\widetilde\F$ as follows: for the generators it is given by
\begin{equation*}
  \begin{aligned}
    \widetilde\Delta{x_1} &= 1\otimes x_1 + x_1 \otimes 1\\
    \widetilde\Delta{x_2} &= x_2\otimes e^{-x_1} + e^{x_1}\otimes x_2\\
    \widetilde\Delta{x_3} &= x_3\otimes e^{-x_1} + e^{x_1}\otimes x_3\\
    \widetilde\Delta{e^{\pm x_1}} &= e^{\pm x_1}\otimes e^{\pm x_1}
  \end{aligned}
\end{equation*}
For arbitrary products of generators, We extend $\widetilde\Delta$ as an algebra homomorphism with respect to $*$. An easy computation shows that $\widetilde\Delta$ is coassociative.

Now we are interested in the algebra $\F = \widetilde\F / I_\varepsilon$, where $I_\varepsilon$ is the two-sided ideal that is generated by the following relations:
\begin{align*}
  &[x_1,x_2] - 2\varepsilon x_2  && e^{\pm x_1}*x_2-e^{\pm2\varepsilon}x_2*e^{\pm x_1}\\ 
  &[x_1,x_3] + 2\varepsilon x_3  && e^{\pm x_1}*x_3-e^{\mp2\varepsilon}x_3*e^{\pm x_1}\\
  &[x_2,x_3] - 2\varepsilon A(\varepsilon^2)\sinh(2x_1) && [x_1,e^{\pm x_1}]
\end{align*}
where $[a,b]:=a*b-b*a$.

 For the coproduct $\Delta$ on $\F$ descending from $\widetilde\Delta$ to be well defined, it is sufficient to prove that $I_\varepsilon$ is a coideal of $\widetilde\F$ (i.e. $\widetilde\Delta(I_\varepsilon) \subset I_\varepsilon\otimes\widetilde\F + \widetilde\F\otimes I_\varepsilon$). This can be verified by a straightforward computation. We only show one case:
\begin{align*}
  &\widetilde\Delta(x_2*x_3-x_3*x_2-2\varepsilon A(\varepsilon^2)\sinh(2x_1))\\[1ex]
  =&\; (e^{x_1}\otimes x_2 + x_2\otimes e^{-x_1})*(e^{x_1}\otimes x_3 + x_3\otimes e^{-x_1})\\
  &\;- (e^{x_1}\otimes x_3 + x_3\otimes e^{-x_1})*(e^{x_1}\otimes x_2 + x_2\otimes e^{-x_1})\\
  &\;- \varepsilon A(\varepsilon^2)(e^{2x_1}\otimes e^{2x_1}-e^{-2x_1}\otimes e^{-2x_1})\\[1ex]
  =&\;e^{2x_1}\otimes (x_2*x_3-x_3*x_2-2\varepsilon A(\varepsilon^2)\sinh(2x_1))\\
  &\;+(x_2*e^{x_1}-e^{-2\varepsilon}e^{x_1}*x_2)\otimes e^{-x_1}*x_3\\
  &\;+e^{-2\varepsilon} e^{x_1}*x_2 \otimes (e^{-x_1}*x_3-e^{2\varepsilon} x_3*e^{-x_1})
\end{align*}
\\
The following Lemma shows that our algebra is a flat deformation of $\R[x_1,x_2,x_3,e^{\pm x_1}]$ :
\begin{lem}
  $\{x_1^{*n_1}*x_2^{*n_2}*x_3^{*n_3}*e^{mx_1}\,|\,(n_i\in\N,m\in\Z\}$ is a basis for $\F$. If $m$ is negative we read $e^{mx_1}$ as $(e^{-x_1})^{-m}$.
\end{lem}
\begin{proof}We check the following facts:
  \begin{itemize}
  \item [i)] Each word can be ordered using the commutation relations. It is clear that this can always be done recursively.
  \item [ii)] The ordered words are linearily independent as they are (commutative) monomials for $\epsilon=0$.
  \end{itemize}
\end{proof}

\begin{lem}
  The coproduct for $\F$ defined above is a non-trivial deformation of the coproduct on $SL(2,\R)^*$ that comes from the group structure.
\end{lem}
\begin{proof}
  For $\varepsilon=0$, our algebra is the algebra of functions on $SL(2,\R)^*$ that are polynomial in the generators. The product is the pointwise commutative product and the coproduct is the natural coproduct comming from the group structure.\\
To see that the coproduct is indeed a non-trivial deformation, we could for example compute $\Delta((x_2)^2)$ ($x_2^{*2}=(x_2)^2$) and we would find that $\Delta(x_2)\Delta(x_2) - \Delta(x_2)*\Delta(x_2)$ does not vanish in second order in $\varepsilon$.
\end{proof}

\section{Bialgebra structure}

\subsection{$U_h(sl(2,\R))$}
So far we constructed a bialgebra $\F = \R\{x_1,x_2,x_3,e^{\pm x_1}\}[[\varepsilon]] / I_\varepsilon$ where $I_\varepsilon$ are the relations given by the commutators of the generators. We would like to write this algebra in different coordinates to see that it is isomorphic to the quantum universal enveloping algebra $U_h(sl(2,\R))$. Therefore we set
\begin{align*}
  x_1 &= \varepsilon z_1\\
  x_2 &= \sqrt{2\varepsilon A(\varepsilon^2)\sinh(h)}\; z_2\\
  x_3 &= \sqrt{2\varepsilon A(\varepsilon^2)\sinh(h)}\;z_3
\end{align*}
and $h=2\varepsilon$. For the new generators $z_i$ the
commutation relations look as follows:
\begin{align*}
  [z_1,z_2] &= \phantom{-}2z_2 & e^{\pm hz_1/2}*x_2 &= e^{\pm h}x_2*e^{\pm hz_1/2}\\
  [z_1,z_3] &= -2z_3           & e^{\pm hz_1/2}*x_3 &= e^{\mp h}x_3*e^{\pm hz_1/2}\\
  [z_2,z_3] &= \frac{\sinh(hz_1)}{\sinh(h)} & [z_1,e^{\pm hz_1/2}]&=0
\end{align*}
and we find the following expressions for the coproduct
\begin{align*}
  \Delta z_1 &= 1\otimes z_1 + z_1\otimes 1\\
  \Delta z_2 &= z_2\otimes e^{-h z_1/2}+e^{h z_1/2}\otimes z_2\\
  \Delta z_3 &= z_3\otimes e^{-h z_1/2}+e^{h z_1/2}\otimes z_3\\
  \Delta e^{\pm hz_1/2} &= e^{\pm hz_1/2}\otimes e^{\pm hz_1/2}
\end{align*}
We see that $z_{1,2,3}$ are the generators of $U_h(sl(2,\R))$.

\subsection{Deformation with 2 Parameters}
We now rewrite our algebra with $x_{1,2,3}$-coordinates defined as follows:
\begin{align*}
  x_1 &= \frac{h}{2}\,\xi_1\\
  x_2 &= \sqrt{A(\varepsilon^2)h\sinh(h)}\,\xi_2\\
  x_3 &= \sqrt{A(\varepsilon^2)h\sinh(h)}\,\xi_3
\end{align*}
We multiply the Poisson structure by $\frac{h}{2}$ and the relations of the algebra become
\begin{align*}
  [\xi_1,\xi_2] &= \phantom{-}2\varepsilon \xi_2 & e^{\pm h\xi_1/2}*\xi_2 &= e^{\pm\varepsilon h}x_2*e^{\pm hz_1/2} \\
  [\xi_1,\xi_3] &= -2\varepsilon \xi_3           & e^{\pm h\xi_1/2}*\xi_3 &= e^{\pm\varepsilon h}x_3*e^{\pm hz_1/2}\\
  [\xi_2,\xi_3] &= \varepsilon\frac{\sinh(h\,\xi_1)}{\sinh(h)} & [\xi_1,e^{\pm h\xi_1/2}] &= 0
\end{align*}
and we find the following expressions for the coproduct:
\begin{align*}
  \Delta \xi_1 &= 1\otimes \xi_1 + \xi_1\otimes 1\\
  \Delta \xi_2 &= \xi_2\otimes e^{-h \xi_1/2}+e^{h \xi_1/2}\otimes \xi_2\\
  \Delta \xi_3 &= \xi_3\otimes e^{-h \xi_1/2}+e^{h \xi_1/2}\otimes \xi_3\\
  \Delta e^{\pm h\xi_1/2} &= e^{\pm h\xi_1/2}\otimes e^{\pm h\xi_1/2}
\end{align*}

\noindent Markus R. Engeli, 
D-MATH, ETH-Zentrum\\
CH-8092 Z\"urich, 
Switzerland\\
{\em e-mail address:} {\tt engeli@math.ethz.ch}

\end{document}